\documentclass[12]{article}

\usepackage{srcltx}
\usepackage{amssymb,amsmath}
\usepackage{color, graphics}
\usepackage[latin1]{inputenc}
\usepackage[dvips]{graphicx}
\usepackage{setspace}
\usepackage[comma,authoryear]{natbib}

 \newtheorem{prop}{Proposition}
 \newtheorem{Lemma}{Lemma}
 \newtheorem{thm}{Theorem}
 \newtheorem{cor}{Corollary}
 \newtheorem{df}{Definition}

 \newcommand{\bt}{\begin{thm}}
  \newcommand{\et}{\end{thm}}
  \newcommand{\bdf}{\begin{df}}
  \newcommand{\edf}{\end{df}}
  \newcommand{\bp}{\begin{prop}}
  \newcommand{\ep}{\end{prop}}
  \newcommand{\bc}{\begin{cor}}
  \newcommand{\ec}{\end{cor}}
  \newcommand{\bl}{\begin{Lemma}}
  \newcommand{\el}{\end{Lemma}}

  \newcommand{\be}{\begin{eqnarray*}}
  \newcommand{\ee}{\end{eqnarray*}}
  \newcommand{\ben}{\begin{eqnarray}}
  \newcommand{\een}{\end{eqnarray}}
  \newcommand{\bit}{\begin{itemize}}
  \newcommand{\eit}{\end{itemize}}
  \newcommand{\st}{\vspace{-0.8cm}\begin{flushright} \mbox{$\diamond$} \end{flushright}}

\newcommand{\Rex}{\mathbb{R}}

 \newcommand{\gmu}{\boldsymbol{\mu}}
\newcommand{\te}{\theta}
\newcommand{\Te}{\Theta}
\newcommand{\teh}{\hat{\theta}}
\newcommand{\Teh}{\hat{\Theta}}
\newcommand{\f}{\phi}

\newcommand{\la}{\lambda}
\newcommand{\al}{\alpha}
\newcommand{\etah}{\hat{\eta}}
\newcommand{\rhoh}{\hat{\rho}}

  \newcommand{\gte}{\boldsymbol{\theta}}
  
  \newcommand{\gfi}{\boldsymbol{\phi}}
  
  \newcommand{\gla}{\boldsymbol{\lambda}}

  \newcommand{\gx}{{\bf x}}
  \newcommand{\gX}{{\bf X}}

  \newcommand{\gG}{{\bf G}}
  \newcommand{\gU}{{\bf U}}
 \newcommand{\gS}{{\bf S}}
  \newcommand{\gs}{{\bf s}}
  \newcommand{\gp}{{\bf p}}
 
  \newcommand{\gV}{{\bf V}}

\newcommand{\lpp}{\ell''(\hat{\te})}
\newcommand{\lppp}{\ell'''(\hat{\te})}

\begin{document}

\begin{center}
{\bf \Large
Some asymptotic results for fiducial and confidence distributions}

\vspace{4mm}

{\Large Piero Veronese \footnotemark and Eugenio Melilli \\

{\small Bocconi University, via R\"ontgen, 1, 20136, Milan, Italy\\
piero.veronese@unibocconi.it, eugenio.melilli@unibocconi.it}}

\end{center}
\footnotetext{Corresponding author}

\bigskip

 \noindent{\bf Abstract}. Under standard regularity assumptions, we provide simple approximations for \emph{specific} classes of fiducial and confidence distributions and discuss their connections with objective Bayesian posteriors. For a real parameter the approximations are accurate at least to order $O(n^{-1})$. For the mean parameter $\gmu=(\mu_1,\ldots, \mu_k)$ of an exponential family, our fiducial distribution is asymptotically normal and invariant to  the importance ordering of the $\mu_i$'s.

\bigskip

 \noindent{\bf Keywords:}  ancillary statistic, confidence curve, coverage probability,  natural exponential family, matching prior, reference prior.



\medskip

\section{Introduction} \label{sec:_intro}

Confidence and fiducial distributions, often confused in the past,  have recently received a renewed attention by statisticians thanks to several contributions which clarify  the concepts within a  purely frequentist setting and  overcome the lack of rigor and completeness typical of the original formulations.
For a wide and comprehensive presentation of the theory of confidence distributions and a rich bibliography we refer the reader to the book by \citet{Schweder:2016} and to the review paper by \citet{Xie:2013}. This latter also highlights the importance of this theory in meta-analysis, see also \citet{liu2015multivariate}.
For what concerns  fiducial distributions Hannig and his coauthors, starting from the original idea of Fisher, have developed in several papers a \emph{generalized fiducial inference} which is suitable for a large range of situations; see \citet{Hannig:2016} for a complete review on the topic and updated references.

Given a random vector $\gS$  (representing the observations or a sufficient statistic) with distribution indexed by  $\boldsymbol{\eta}=(\te,\gla)$, where $\te$ is the real parameter of interest, a  \emph{confidence distribution} (CD) for $\te$ is a function $C$ of $\gS$ and $\te$ such that: i)  $C(\gs,\cdot)$ is a distribution function on $\Rex$ for any fixed realization $\gs$ of $\gS$ and ii) $C(\gS,\theta)$ has a uniform distribution on $(0,1)$,  whatever the true value of $\boldsymbol{\eta}$.
The second condition is crucial because it implies that the coverage of the intervals derived from $C$ is exact.
 If it  is satisfied only for the sample size tending to infinity, $C$ is an \emph{asymptotic} CD and  the coverage is correct only approximately.
Given a CD, it is possible to define the  \emph{confidence curve} $cc_\gs(\theta)=|1-2C(\gs,\theta)|$, which displays the confidence intervals induced by $C$ for all levels, see \citet[Sec. 1.6]{Schweder:2016}.

A \emph{fiducial distribution} (FD) for a parameter $\gte$ has been obtained by several authors starting from a data-generating equation $\gS=\gG(\gU,\gte)$, with $\gU$  a random vector with known distribution, which allows to transfer randomness  from $\gS$  to $\gte$. In particular \citet[2016]{Hannig:2009} derives an explicit expression for the density of a FD which coincides with that originally proposed by \citet{Fisher:1930}, namely $h_s (\theta)= |\partial F_\theta(s)/ \partial \theta|$, when  both $\te$ and $S$ are real and $G(U, \te)=F_{\te}^{-1}(U)$, with $F_{\te}$ distribution function of $S$ and $U$ uniform in $(0,1)$.

In this paper we consider the specific definition of FD given in \citet{Veronese:2015}, recalled in Section \ref{sec:multiparameter}, which for a real parameter and a continuous $S$ again simplifies to the Fisher's formula. In particular, we assume that the FD function is
\ben \label{h}
H_s(\theta)=1-F_\theta (s) = \mbox{Pr}_{\te}(S>s)
\een
 with  $F_\theta (s)$ decreasing and differentiable  in $\te$ and with limits $0$ and $1$ when $\te$ tends to the boundaries of its parameter space. This conditions are always  true, for example, if $F_{\te}$ belongs to a regular real natural exponential family (NEF).
This FD  is also a CD (asymptotically in the discrete case).
For the multi-parameter case a peculiar aspect of our FD is its dependence on the inferential importance ordering of the parameters, similarly to what happens for the objective Bayesian posterior obtained from a reference prior.
The connections between our definition and Hannig's setup are discussed in \citet{Veronese:2015}.

In Section \ref{sec:first-order}, extending a result proved in \citet{Veronese:2014} for a NEF, we give a second order asymptotic expansion of our FD/CD in the real parameter case based only on the maximum likelihood estimator (MLE). This expansion does not require any other regularity conditions  than the standard ones usually assumed in maximum likelihood asymptotic theory. Furthermore, we show that it coincides with the expansion of the Bayesian posterior induced by the Jeffrey prior. This fact establishes a connection with objective Bayesian inference, whose aim is to produce posterior distributions free of any subjective prior information.
In Section \ref{sec:p-star},  starting from the well known $p^*$\emph{-formula}  of \citet[1983]{Barndorff:1980}, we propose and discuss a FD/CD which, using an ancillary statistic  in addition to the MLE, has good asymptotic behavior. Higher order asymptotics for generalized fiducial distributions have been discussed, at our knowledge, only in the unpublished paper \citet{Majunder:2016}. However, its  focus is different being devoted to identify data generating equation with desirable properties.
In Section \ref{sec:multiparameter} we consider a NEF with a multidimensional parameter and show  that, without any further regularity conditions, the asymptotic FD of the mean parameter is normal, it does no longer depend on the inferential ordering of the parameters  and coincides with the corresponding asymptotic Bayesian  posterior.
Some examples illustrate the good properties and performances of the various proposed FD/CD with emphases on coverage and expected length of confidence intervals. Finally, the Appendix  includes the proofs of all the theorems and propositions stated in the paper.

\section{Asymptotics for fiducial and confidence distributions: the real parameter case} \label{sec:asymptotic}

\subsection{An expansion with error of order $O(n^{-1})$} \label{sec:first-order}

In \citet{Veronese:2014} an Edgeworth expansion with an error of order $O(n^{-1})$  of  the FD/CD  for the mean parameter of a real NEF was derived.
Here we generalize this result to an arbitrary regular model.

 Let $\gX=(\gX_1, \ldots, \gX_n)$ be an i.i.d. sample of size $n$ from a density (with respect to the Lebesgue measure) parameterized by $\te$ belonging  to an open set $\Te \subseteq \Rex$. Let $\teh$ be the MLE of $\te$ based on $\gX$ and denote by $p_\te(\teh)$ its density.
Let  $\ell(\te)=n^{-1}\log p_\te(\teh)$ and let  $\lpp$ and $\lppp$ be the second and the third derivative  of $\ell(\te)$ with respect to $\te$, evaluated in $\teh$. Then the expected and observed Fisher information of $\teh$ (per unit) are $I(\te)=-n^{-1}E_{\te}(\partial^2 \log p_\te(\teh)/\partial \te^2)$ and $-\lpp$, respectively. Let $b=b(\teh)=-1/\lpp$. Consider now $Z=\sqrt{n/ b}\,(\te-\hat{\te})$, which is an approximate standardized version of  $\te$ in the FD/CD-setup, and let $H_{n,\teh}(z)$ be its FD/CD derived from  the sampling distribution of $\hat{\te}$.
If $\teh$ is sufficient, $H_{n,\teh}(z)$ is exact, otherwise it is a natural approximation of the exact one, see e.g. \citet{Schweder:2016}.
To prove our result we resort to  the expansion of the frequentist probability $\mbox{Pr}_{\te}(Z\leq z)$ provided in  \citet{Datta-Ghosh:1995} or in \citet{mukerjee:1997}. Thus we need the regularity assumptions
used in these papers, see also \citet[Ch. 8]{Ghosh:1994} and \citet{Bickel:1990} for a precise statement.
Notice that the conditions required for the frequentist expansion of the distribution of the MLE are rarely reported in a rigorous way in books and papers. However, what is important here is that, in order to prove our result, we do not need any further assumption  and this fact allows an immediate and fair  comparison between MLE- and FD/CD-asymptotic theory.
\bt \label{t:asymptotic}Let $\gX$ be an i.i.d. sample of size $n$ from a  density $p_{\te}$, $\te \in \Te \subseteq \Rex$. Then, under the regularity assumptions cited above,  the distribution function $H_{n,\teh}(z)$ of the FD/CD for $Z=\sqrt{n/ b}(\te-\hat{\te})$ has the expansion
\ben \label{expansion teta}
H_{n,\teh}(z)=\Phi(z)-\phi(z)\left[ \frac{1}{6}b^{3/2}\lppp(z^2-1)\right]n^{-1/2} + O(n^{-1}).
\een
\et
If $p_{\te}$ also satisfies the conditions for the expansion of a Bayesian posterior, see e.g.  \citet[Theorem 2.1, with $K=1$]{Johnson:1970},  we have the following

\bc \label{c:Jeffreys-exp}
 In the same setting of Theorem \ref{t:asymptotic}, let $\pi^J(\te)\propto I(\te)^{1/2}$ be the Jeffreys prior for $\te$. If $\pi^J$ is improper, assume that  there exists an $n_0 \geq 1$ such that the posterior distribution $\pi^J(\te|z)$ of $\te$ is proper for $n\geq n_0$, almost surely for all $\te$. Then the expansion of $\pi^J(\te|z)$ coincides with that of  $H_{n,\teh}(z)$ given in \eqref{expansion teta}.
\ec

Theorem \ref{t:asymptotic} and Corollary \ref{c:Jeffreys-exp} confirm the idea that the Jeffreys posterior is really free of any subjective prior information. Furthermore, they naturally establish a connection between FD/CD-theory and \emph{matching priors}, i.e. priors that ensure approximate frequentist validity of posterior credible sets. More precisely,  a prior $\pi$ for which
$\mbox{Pr}_{\te}(\te\leq q_{1-\al}(\gX,\pi))=1-\al +o(n^{-r/2})$, where  $q_{1-\al}(\gX,\pi)$ denotes the $(1-\al)$th posterior quantile of $\te$, is called a matching prior of order $r$, see \citet{Datta-Muk:2004} for a  general review and references.
For a regular model indexed by a real parameter  it is well known, see \citet[Theorem 2.5.1]{Datta-Muk:2004}, that the Jeffreys prior $\pi^J$ is the unique first order matching prior and is also a second order matching prior if and only if the model satisfies the following  condition:
\ben \label{second-order-cond}
I(\te)^{-3/2}E_\te[(\partial\ell(\te)/\partial \te)^3] \quad \mbox{ is a constant free of $\te$.}
\een
\citet{Veronese:2014} study the existence of a prior (named \emph{fiducial prior}) which induces  a Bayesian posterior coinciding with the FD, extending a result given by \citet{Lindley:1958} for a continuous univariate sufficient statistic, see also \citet{Taraldsen:2015} for a generalization to multivariate group models.
 Because a FD/CD realizes the exact matching, we immediately have  the following
\bc \label{cor: fiducial prior}
 If a fiducial prior $\pi^F$ exists, then it coincides with the Jeffreys prior $\pi^J$. Furthermore, the condition \eqref{second-order-cond} is necessary for the  existence of $\pi^F$.
 \ec
 Notice that for a model belonging to a NEF, with  mean parameter $\mu$ and variance function $V(\mu)$, condition \eqref{second-order-cond} becomes:
 \lq\lq$2V'(\mu)V(\mu)^{-1/2}$ is constant\rq\rq{}. The solution of this differential equation is  $V(\mu)=(c_1 \mu+c_2)^2$, i.e. a fiducial prior for a parameter of a NEF may exist only if its variance function is quadratic. This result was found for the first time in \citet{Veronese:2014}, using a totally different approach.

\begin{figure}
\vspace{-1.5cm}
\centering
\includegraphics[width=.49\textwidth,height=5cm]{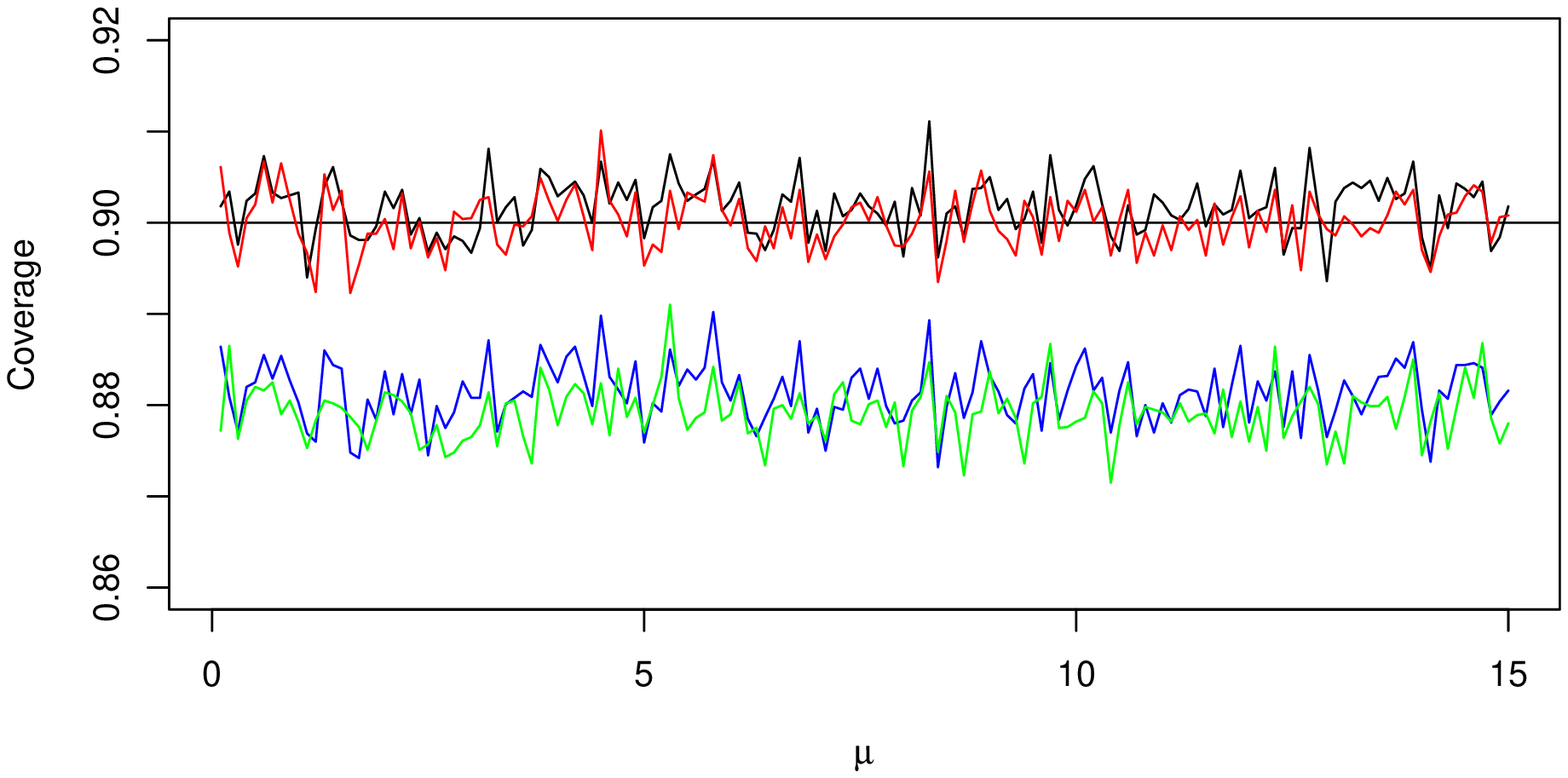}\hfil
\includegraphics[width=.49\textwidth,height=5cm]{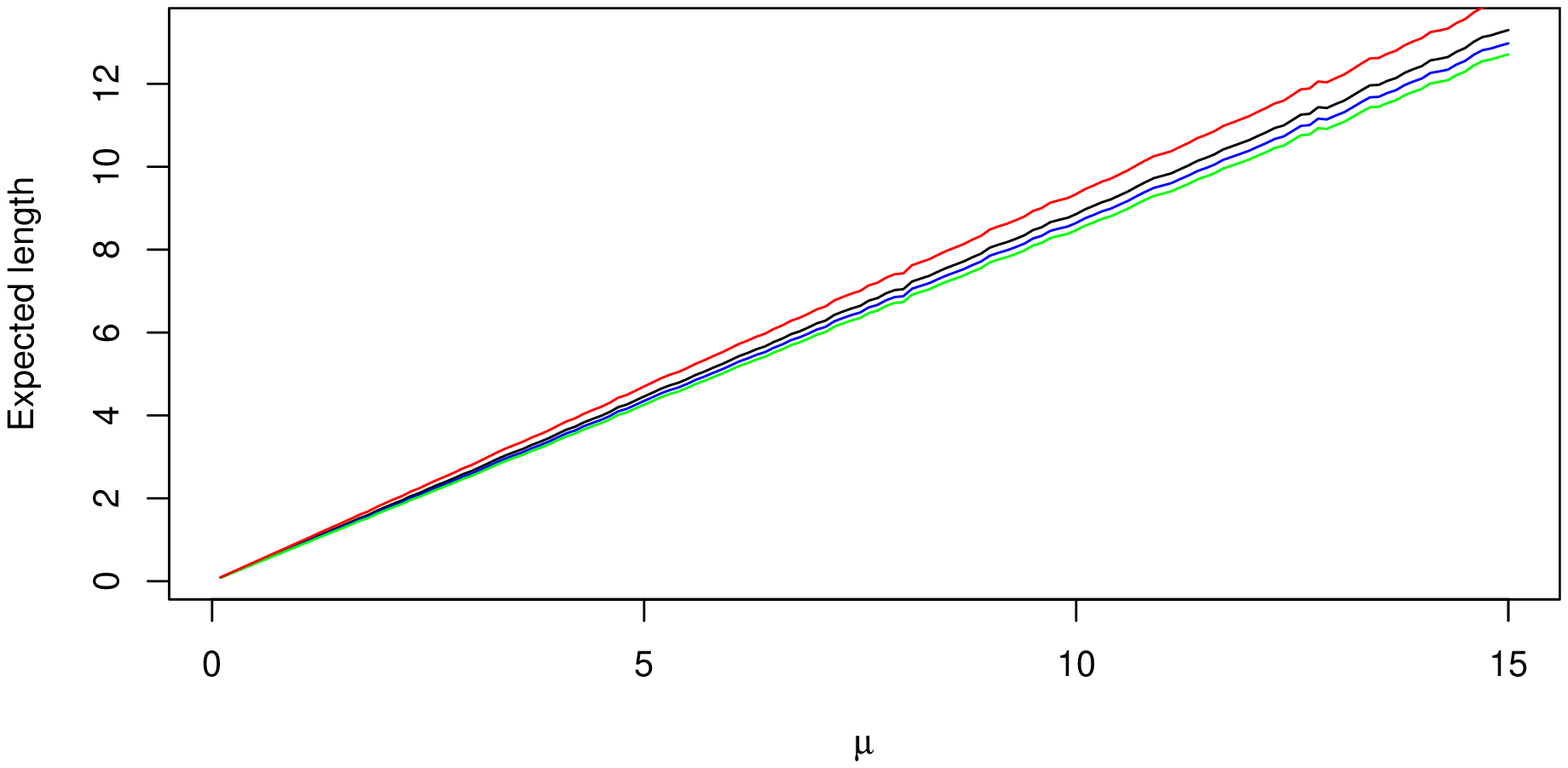}
\parbox{14cm}{
\vspace{-0.5cm}\caption{\small{Coverages and expected lengths for the 90\%  intervals with $n=15$ based on: exact FD/CD (red), Normal approximation (green), expansions of FD/CD (black) and of MLE (blue). }}\label{fig:exp-neg_intervals}}
\end{figure}

\emph{Example 1} (Exponential distribution). 
Let $(X_1,\dots,X_n)$ be an i.i.d. sample from an exponential distribution with mean $\mu$. The MLE  $\hat{\mu}$ of $\mu$ is the sample mean and $b=b(\hat{\mu})=\hat{\mu}^2$. Then the expansions of the FD/CD for $Z=(\sqrt{n}/\hat{\mu})(\mu-\hat{\mu})$ in \eqref{expansion teta} and that of the standardized MLE $W=-Z=(\sqrt{n}/\hat{\mu})(\hat{\mu}-\mu)$, see (A.3), are respectively
\ben \label{exp-neg-approx_1}
\Phi(z)-\phi(z)\left[ 2(z^2-1)/3\right]n^{-1/2}+ O(n^{-1}),\\
\Phi(w)-\phi(w)\left[ -(2w^2+1)/3\right]n^{-1/2}+ O(n^{-1})  \label{exp-neg-approx_2}.
\een
It follows that the confidence intervals obtained from \eqref{exp-neg-approx_1} and \eqref{exp-neg-approx_2} are different, contrary to what happens for those based only on  the normal approximation.  Their coverages and expected lengths are reported in Figure~\ref{fig:exp-neg_intervals} for a sample of size $n=15$ and confidence level 0.9. Notice that  the coverage of the FD/CD-intervals is much closer to the nominal level than that of the intervals based on the MLE, while the expected lengths are quite similar. For the sake of comparison Figure \ref{fig:exp-neg_intervals} reports also the coverage and the expected length of the intervals based on the exact FD/CD, which is an  inverse-gamma($n,n\hat{\mu})$, see \citet[Tab.1]{Veronese:2014}. The latter intervals are clearly exact, but wider.
Finally, by  Corollary \ref{c:Jeffreys-exp},  the expansion \eqref{exp-neg-approx_1} coincides with that of the Jeffreys posterior. It is easy to verify, according to Corollary \ref{cor: fiducial prior}, that the fiducial prior exists and coincides with  $\pi^J(\mu)\propto 1/\mu$.
\st

\begin{figure}
\vspace{-1.0cm}
\begin{center}
\includegraphics[width=7cm,height=3.5cm]{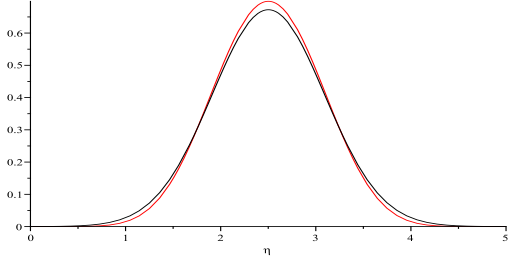}\\
\parbox{13cm}{
\vspace{-0.3cm}
\caption{\small{Approximate (black) and exact(red) fiducial densities for $\eta$ in the Fisher's gamma hyperbola for a sample size $n=5$, $s_1=17.321$, $s_2=0.116$. \label{fig:nilo}}}}
\end{center}
\vspace{-0.5cm}
\end{figure}

\emph{Example 2} (Fisher's gamma hyperbola-Nile problem). \citet[Sec.VI.9]{Fisher:1973} considers a sample of size $n$ from a curved exponential family obtained by two independent gamma distributions with means constrained on an hyperbole.  Following \citet{Efron-Hinkley:1978}, we directly start with the sufficient statistic $\gS=(S_1,S_2)$,  with $S_1$ and $S_2$ distributed according to  ga$(n,e^{-\eta})$ and ga$(n,e^{\eta})$, respectively. Here ga$(\alpha,\beta)$ denotes a gamma distribution with shape parameter $\alpha$ and mean $\alpha/\beta$. It follows that the likelihood of the model is
$L_{\eta}(\gs)=\exp\{-e^{-\eta}s_1-e^{\eta}s_2\}$ and that  the
 MLE of $\eta$ is $\etah=(1/2) \log(S_1/S_2)$.
 Even if an exact inference on $\eta$ cannot be performed using only $\etah$, the minimal sufficient statistic $\gS$ is indeed bivariate, an asymptotic FD/CD for $\eta$, based on $\etah$, can be easily obtained  from Theorem \ref{t:asymptotic}.
Since $\ell'''(\hat{\eta})=0$, it follows from \eqref{expansion teta} that, in this case, the normal distribution N$(\etah, b/n)$, with $b=-1/\ell''(\hat{\eta})=n/(2\sqrt{s_1s_2})$, is an approximate FD/CD of $\eta$ with error of order $O(n^{-1})$. Figure \ref{fig:nilo} reports the plot of its density compared with the exact FD/CD based on $\gS$, which will be derived in the next section. It shows the goodness of the approximation  even  for a very small sample size ($n=5$ in the plot).
Finally, it is easy to check that $E_\eta[(\partial \ell(\eta)/\partial \eta)^3]=0$, thus condition \eqref{second-order-cond} holds and, by Corollary \ref{cor: fiducial prior}, a fiducial prior might  exist for this model. Indeed it exists and we will find it in the next section.  \st

Another criterium to define  matching priors studied in the Bayesian literature is based directly on the distribution functions, see \citet[Sec.\hspace{0.1 cm}3.2]{Datta-Muk:2004}.  Because $H_{n,\teh}(z)$ is stochastic in a frequentist setup, as it occurs for a posterior distribution, we can consider the matching between $E_\te(H_{n,\teh}(z))=E_\te\left(\mbox{Pr}_{\teh}\left \{\sqrt{n/b}(\te -\teh)\leq z \right \}\right)$ and $\mbox{Pr}_{\te}\left \{\sqrt{n/b}(\te -\teh)\leq z \right \}$. Clearly quantiles and distribution functions are  strongly connected and thus it  is not  surprising that the conditions for the existence of matching priors in the two criteria are related. Indeed, the first order matching conditions are the same, while this is not true for the second order ones.
Notice that the matching in terms of quantiles is obtained using the quantity $Z$ which can be seen as an approximate pivotal quantity. This is meaningful in an asymptotic setting, but it is not appropriate for small sample sizes. In this case,  the FD/CD realizes an exact matching if we replace  $Z$ with the pivotal quantity given by the distribution function of $\teh$, namely $F_\te(\teh)$. Indeed, we have
 \be
& & E_\te \left(\mbox{Pr}_{\teh}\{F_\te(\teh)\leq z\} \right)= E_\te \left(\mbox{Pr}_{\teh}\{1-H_{\teh}(\te)\leq z\} \right) =\\
& & E_\te \left(\mbox{Pr}_{\teh}\{H_{\teh}(\te)\geq 1-z\} \right)
=1- E_\te \left(\mbox{Pr}_{\teh}\{H_{\teh}(\te)\leq 1-z\} \right)=1-E_\te (1-z)=z,
\ee
and because $ \mbox{Pr}_{\te}\{F_\te(\teh)\leq z\} =z$, the exact matching for distribution functions holds.
However, an exact FD/CD does not always exist and thus it is natural to look for approximations which have nice asymptotic properties.  Furthermore, in a multiparameter case quantiles are not well defined and thus the study of the frequentist properties of a multivariate FD/CD can be conducted along the lines developed for matching distribution functions.

\subsection{An approximation based on the Barndorff-Nielsen $p^*$-formula} \label{sec:p-star}

Consider a sample $\gX$  whose distribution depends on a real parameter $\te$. In the previous section we have obtained an approximate FD/CD for $\te$ starting from the distribution of the MLE  $\teh$. However,  if $\teh$ is not sufficient,  the approximation of the FD/CD can be improved adding the remaining information included in the sample. This can be done resorting  to the \lq\lq conditionality resolution\rq \rq{} of the statistical model, i.e. the construction of an ancillary statistic $\textbf{A}$ and of an approximate conditional distribution of $\teh$ given $\textbf{A}=\textbf{a}$. We refer to \citet[1983]{Barndorff:1980} for a detailed discussion on the topic and recall here only some useful facts.
His well known  approximate distribution of $\teh$ given $\textbf{A}=\textbf{a}$  is
\ben \label{pstar}
p^*_\te(\teh|\textbf{a})=c(\textbf{a},\te)|j(\teh)|^{1/2} L(\te;\gx)/L(\teh;\gx),
\een
where  $L(\te;\gx)$ is the likelihood function,  $ j(\teh)$ is the observed Fisher information and $c(\textbf{a},\te)$ is the normalizing constant which  does not depend on $\te$ in many important cases.
Formula \eqref{pstar} is quite simple,  is generally accurate to order $O(n^{-1})$, or even $O(n^{-3/2})$, and exact in specific cases.  Here the term approximation refers to one of the two following situations: i) there exists an ancillary statistic $\textbf{A}$, but it is not possible to construct the exact conditional distribution of $\teh$ given $\textbf{A}=\textbf{a}$; ii) an exact ancillary statistic does not exist and an approximate one is used.
It is worth to remark that formula $p^*$ is  invariant to reparameterizations and  is exact for transformation models. Furthermore, under repeated sampling from a real NEF, where no conditioning is involved, $p^*$  is often of order $O(n^{-3/2})$ and  is exact for normal  (known variance), gamma (known shape) and inverse-gaussian (known shape) distributions.

If $F^*_{\te}(\teh|\textbf{a})$ denotes the distribution function corresponding to $p^*_\te(\teh|\textbf{a})$ and  satisfies the conditions reported after \eqref{h}, we can derive an approximate FD/CD for $\te$ as $h^*_{\gx}(\te)=|\partial F^*_{\te}(\teh|\textbf{a})/\partial \te| $.
 This construction of a FD/CD  is not different in essence from the widespread procedure used to derive a Bayesian posterior starting from an approximate (e.g. profile,  pseudo  or  composite) likelihood. A similar approach based on approximate likelihood is used also by \citet{Schweder:2016} to construct a CD.

 The next result concerning a real NEF is useful when the exact distribution of $\teh$  is difficult to obtain.

\bp \label{prop_pstar}
If $\teh$ is the MLE of $\te$ based on an i.i.d. sample from a real regular NEF, with density $p_\te(x)= \exp\{\te x -M(\te)\}$, then
$h^*_{\teh}(\te)=|\partial F^*_{\te}(\teh)/\partial \te| $ is an exact FD/CD for $\te$ based on  $p^*_{\te}(\teh)$. It is an approximate FD/CD based on the whole sample and its order of approximation depends on that of $p^*_{\te}(\teh)$.
\ep

The following examples, concerning  curved exponential families, i.e.  NEFs in which a constraint on the natural parameter space is imposed, illustrate another typical case in which formula \eqref{pstar} can be fruitfully applied to construct a FD/CD.

\emph{Example 2 ctd.}  As previously observed, the MLE $\etah$ is not sufficient and thus the exact FD/CD can be obtained starting from the conditional distribution of $\etah$ given the ancillary statistic $A=\sqrt{S_1 S_2}/n$, proposed by \citet[Sec. VI.10-11]{Fisher:1973}. After some  calculations, one obtains
\ben \label{pstar-iperbole}
p_{\eta}(\etah|a)=\exp\{-2na \cosh(\etah-\eta)\}/(2 K_0(2na)),
\een
where $K_0(w)=\int_{0}^{\infty}\exp\{-w \cosh(z)\}dz$ is the modified Bessel function of the second order evaluated in $(0,w)$.
 As observed by  \citet{Efron-Hinkley:1978},  it is easy to see from \eqref{pstar-iperbole} that this example involves  a  translation (and thus a transformation) model, so that $p_{\eta}(\etah|a)=p^*_{\eta}(\etah|a)$.
 Thus the exact FD for $\eta$ is $h_{\etah,a}(\eta)=-\partial F^*_{\eta}(\etah|a)/\partial \eta $ and, because $\eta$ is a location parameter, it equals the posterior obtained from the Jeffreys prior $\pi^J(\eta) \propto 1$, see \citet[Prop.8]{Veronese:2015}. The nature of the parameter $\eta$ also implies that inferences based on MLE and $h_{\etah,a}(\eta)$ coincide. \st

\emph{Example 3} (Bivariate normal model).
Consider an i.i.d. sample $(X_i, Y_i)$, $i=1, \ldots, n$,  from  a bivariate normal distribution  with expectations 0, variances 1 and correlation coefficient  $\rho$. This is a simple curved exponential model with sufficient statistics $S_1=\sum_{i=1}^n ( X_i^2+Y_i^2)/2$ and $S_2=\sum_{i=1}^n X_iY_i$, but the inference on $\rho$ is a challenging problem as shown in \citet{Fosdick-Raftery:2012} and \citet{Fosdick-Perlman:2016}. Both  \citet{Efron-Hinkley:1978} and \citet{Barndorff:1980} use this example to illustrate the construction of an approximate ancillary statistic in a conditional inference setting. Their proposals essentially coincide and lead to consider the \lq\lq affine\rq\rq{} ancillary  $A=(S_1-n)/\sqrt{n(1+\rhoh^2)}$, where $\rhoh$ is the MLE of $\rho$.

\begin{figure}
\vspace{-1.5cm}
\centering
\includegraphics[width=.49\textwidth,height=5cm]{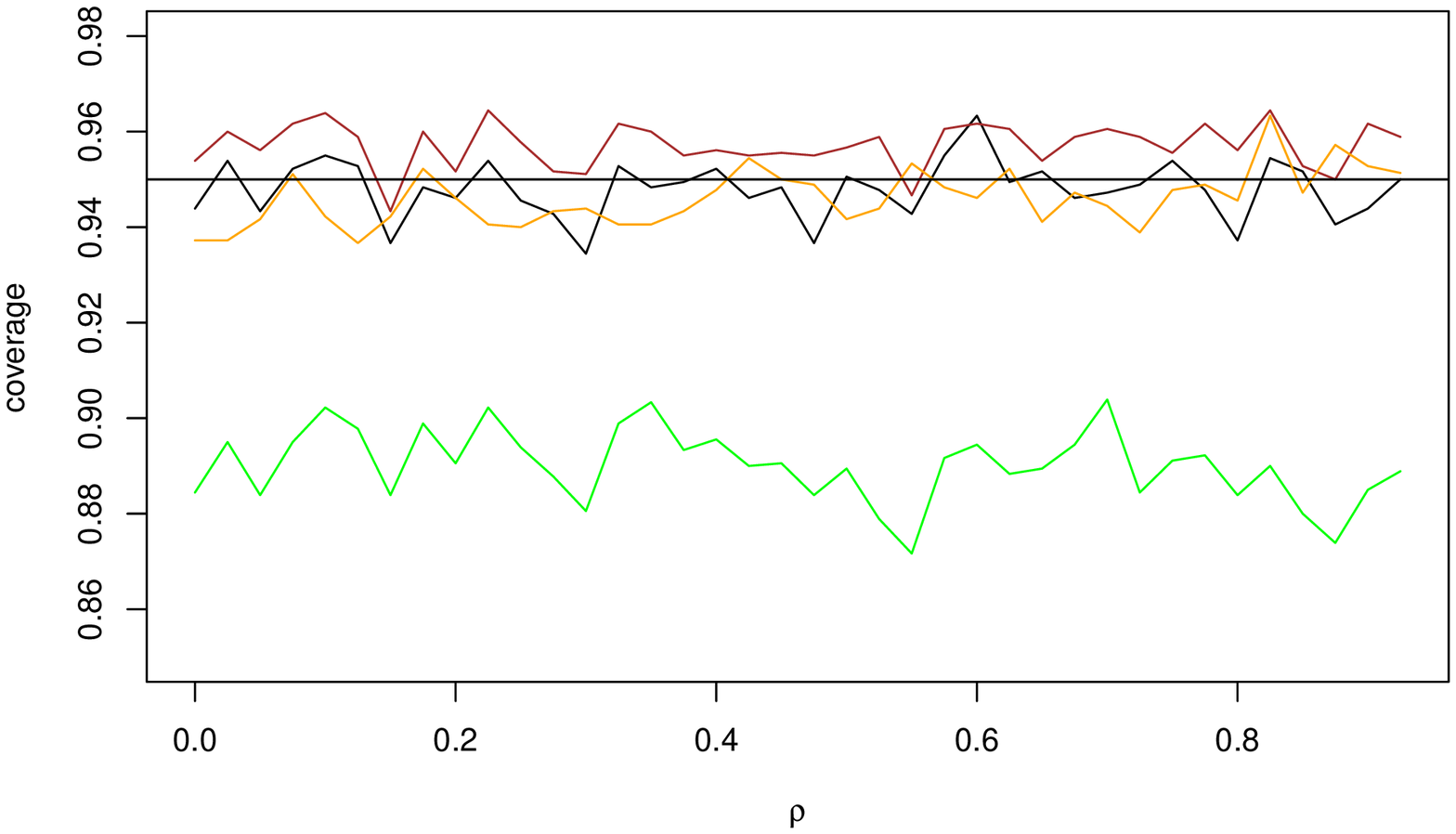}\hfil
\includegraphics[width=.49\textwidth,height=5cm]{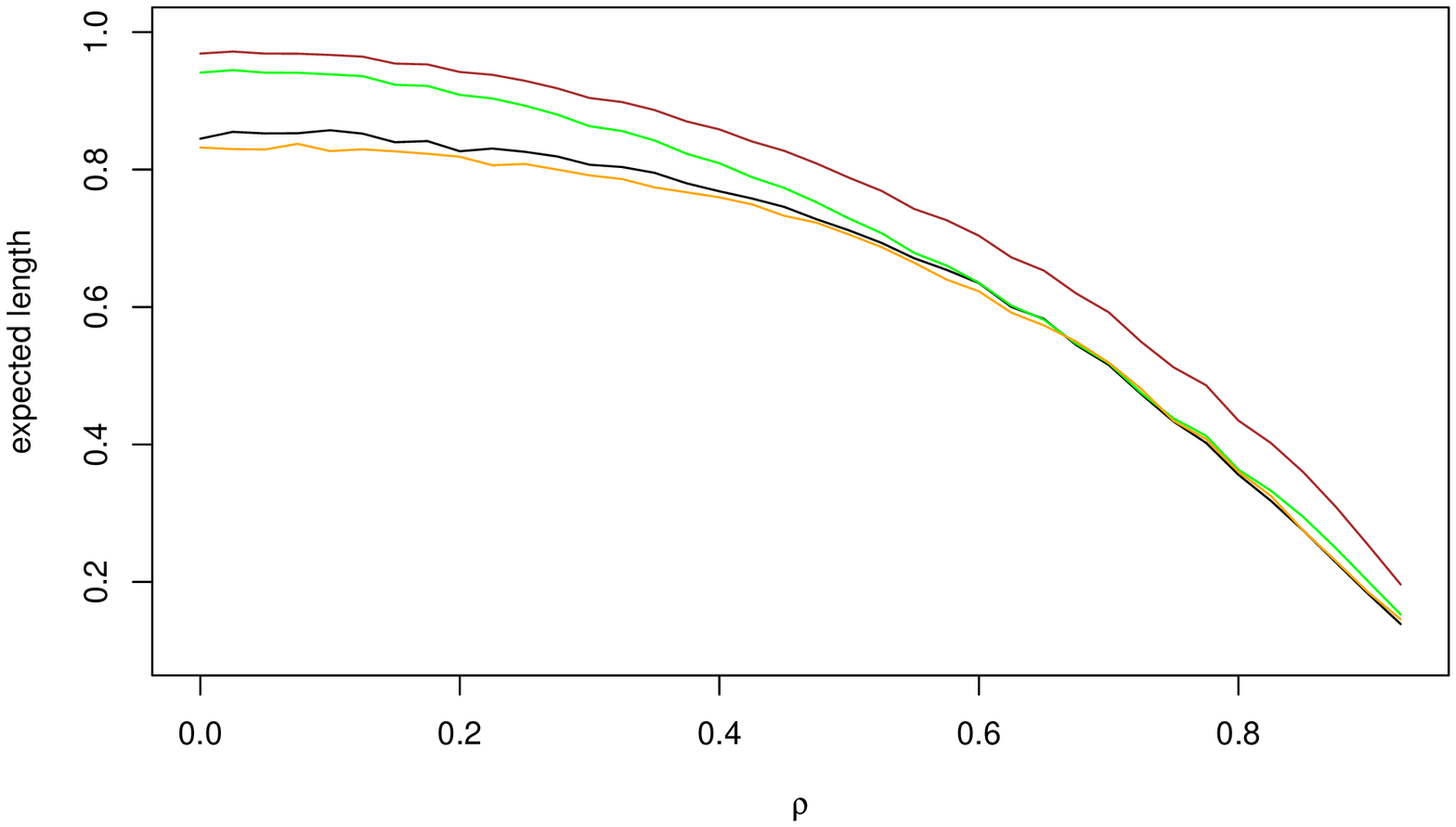}
\parbox{14cm}{
\vspace{-0.5cm} \caption{\small{Coverages and expected lengths of the 95\% intervals with $n=15$ based on: $h^*$ (black), $h^{rstab}$ (brown), $\pi^J$ (orange) and $h^r$ (green).}}\label{fig:rho_intervals}}
\end{figure}

To discuss the performance of $h^*$ obtained starting from $p^*$,  we compare it with other possible asymptotic FDs and with the Bayesian posterior obtained from the Jeffreys prior $\pi^J(\rho) \propto (\rho^2+1)^{1/2}/(1-\rho^2)$.
In particular, we consider the following FDs: $h^r$ and  $h^{rstab}$ obtained from the sample  correlation coefficient $r$ and its stabilizing transformation which, as well known, improves the inferential performance of $r$, see also \citet[pag. 209 and 224]{Schweder:2016}; $h^0$ and $h^1$ obtained considering the first one or the first two terms of \eqref{expansion teta}, respectively. We assume a sample size $n=15$ because a larger value of $n$, e.g. $50$, produces essentially the same (good) results for all choices.
\begin{figure}
\vspace{-1.5cm}
\centering
\includegraphics[width=.49\textwidth,height=4cm]{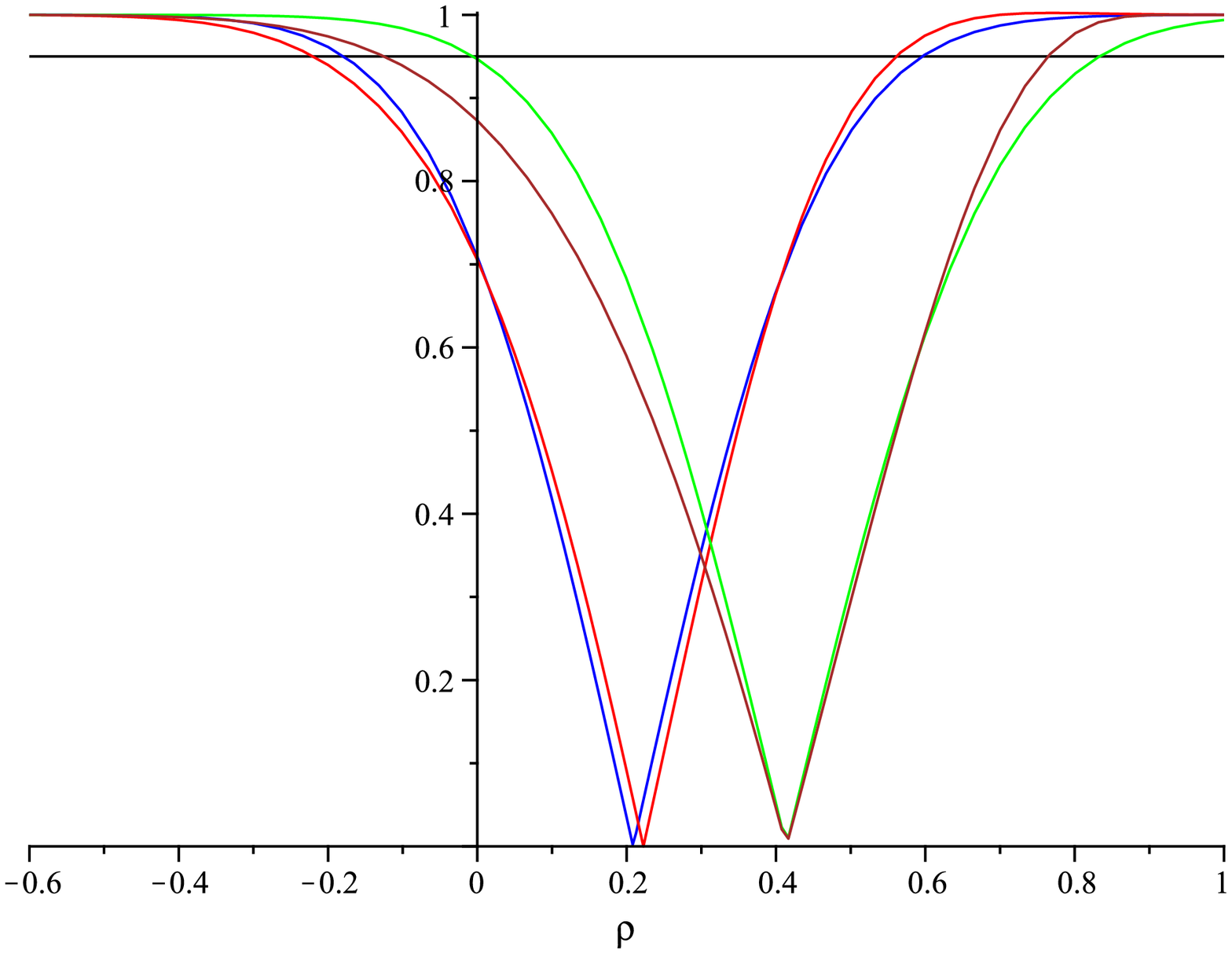}\hfil
\includegraphics[width=.49\textwidth,height=4cm]{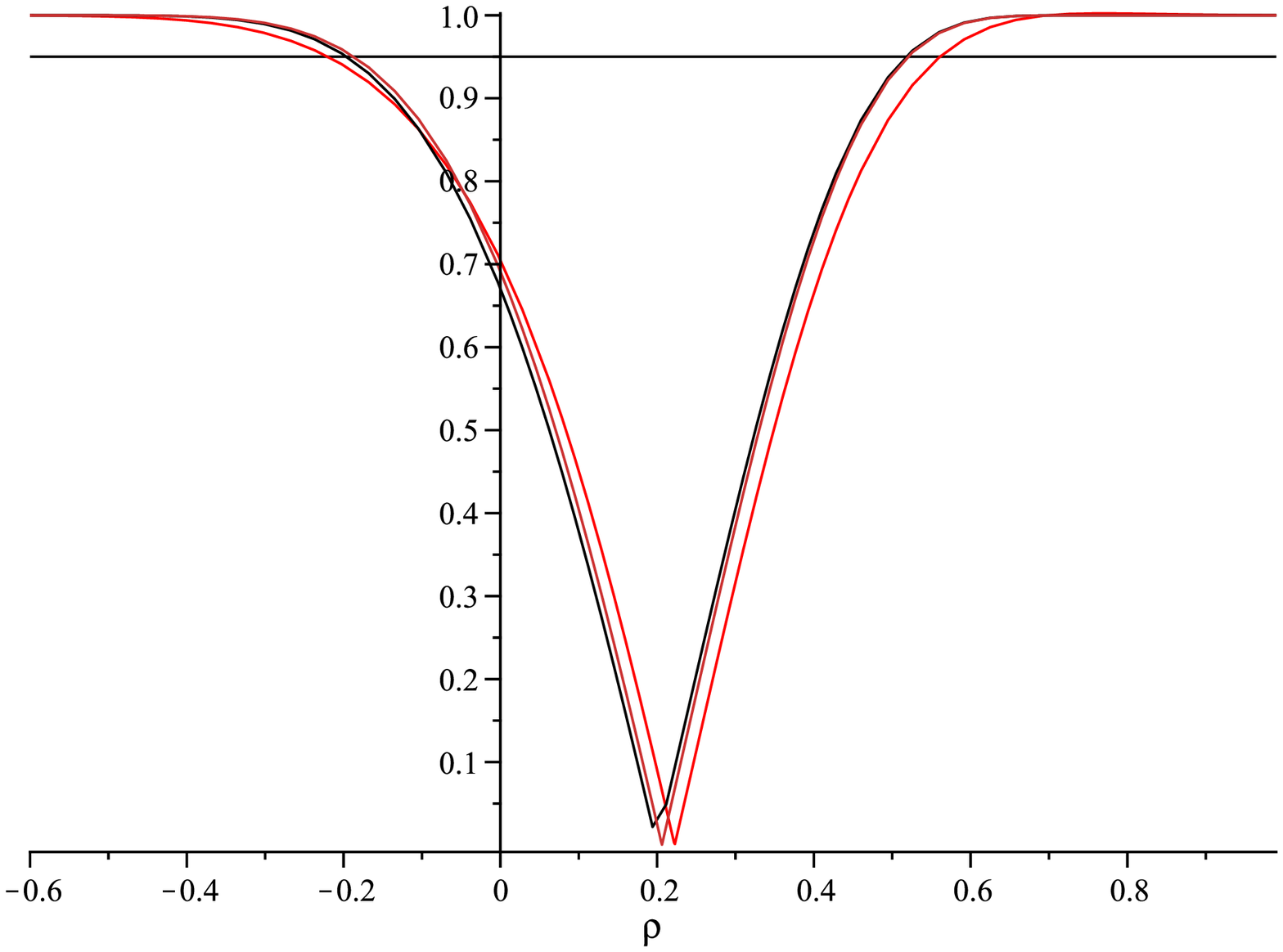}
\parbox{14cm}{
\caption{\small{Confidence curves for a sample size $n=15$, generated from $\rho=0.3$ with $s_1=19.248$ and $s_2=4.827$, $r=0.414$ and $\hat{\rho}=0.209$. Left graph: $cc^r$ (green), $cc^{rstab}$ (brown), $cc^0$ (blue) and $cc^1$ (red). Right graph: $cc^1$ (red), $cc^*$ (black), $cc^J$ (orange). The horizontal line identifies the 95\% confidence intervals.}}\label{fig:rho}}
\end{figure}
The left graph of Figure  \ref{fig:rho} reports an example of the confidence curves $cc^{r}$, $cc^{rstab}$, $cc^0$ and $cc^1$ corresponding to the previous FDs. The curves present different behaviors because they are based on the two estimators $r$ and $\rhoh$ of $\rho$, which assume quite different values in the sample. The right graph compares $cc^1$ with $cc^*$ and $cc^J$ obtained from $h^*$ and Jeffreys posterior, respectively. As expected, the last two curves, both based on the sufficient statistics $S_1$ and $S_2$, are  very  similar and induce    confidence  intervals narrower than those induced by $cc^1$.
To better appreciate the good behavior of $h^*$,  we compare the corresponding coverage and expected length with those of $h^r$, $h^{rstab}$ and the Jeffreys posterior. Figure \ref{fig:rho_intervals} confirms  the very bad inferential performance of $h^r$. The intervals corresponding to $h^*$ have the coverage closest to the nominal one, while those obtained by $h^{rstab}$ present an over-coverage. However, these latter intervals have a uniformly larger expected length. Finally, Bayesian intervals show an intermediate behavior in terms of both coverage and expected length. The same example is discussed by \citet{Majunder:2016}, but they have a different aim and consider different FDs.

\section{Asymptotics for fiducial distributions: the multidimensional parameter case} \label{sec:multiparameter}

For a parameter $\gte$ in $\Rex^d$, inspired by the step-by-step procedure proposed by \citet{Fisher:1973},   \citet{Veronese:2015} give a simple and quite general  definition of FD,  which we summarize here. We refer to the latter paper for details, examples, relationships with objective Bayesian analysis performed using reference priors and a comparison with Hannig's fiducial approach. Notice that for a multidimensional parameter there is not a unique definition of CD, see \citet[Ch.9]{Schweder:2016}, so that in the following we refer only to  FDs.

Given  a random vector $\gS$, representing the sample or a sufficient statistic,  with dimension $m \geq d$ and density $p_{\gte}$, consider the partition $\gS=(\gS_{[d]}, \gS_{-[d]})$, where $\gS_{[d]}=(S_1, \dots, S_d)$ and  $\gS_{-[d]}= (S_{d+1}, \dots, S_m)$,  and suppose that $\gS_{-[d]}$ is ancillary for $\gte$. Clearly, if $d=m$, $\gS_{-[d]}$ disappears. Thus, the density $p_{\gte}$ of $\gS$ can be written as $p_{\gte}(\gs_{[d]}|\gs_{-[d]})p(\gs_{-[d]})$ and the information on $\gte$ provided by the whole sample  is included in the conditional distribution of $\gS_{[d]}$ given $\gS_{-[d]}$.
Assume now that there exists a one-to-one smooth reparameterization from $\gte$ to $\gfi =(\f_1, \ldots, \f_d)$, with the $\f_i$'s ordered with respect to  their inferential importance, such that
\ben \label{general-conditional}
p_{\gfi}(\gs_{[d]}|\gs_{-[d]})=\prod_{k=1}^d p_{\f_{d-k+1}}(s_k|\gs_{[k-1]}, \gs_{-[d]};\gfi_{[d-k]}),
\een
with obvious meaning for $\gs_{[0]}$ and $\gfi_{[0]}$. If,  for each $k$, the one-dimensional conditional distribution function of $S_k$ is monotone and differentiable  in $\f_k$ and has limits $0$ and $1$ when $\f_k$ tends to the boundaries of its parameter space (this is always  true, for example, if this distribution belongs to a regular real  NEF), it is possible to define the joint fiducial density of $\gfi$  as
\ben \label{fid-general}
h_{\gs}(\gfi)=\prod_{k=1}^d h_{\gs_{[k]}, \gs_{-[d]}}(\f_{d-k+1}|\gfi_{[d-k]}),
\een
 where
\ben \label{fid-general-2}
h_{\gs_{[k]}, \gs_{-[d]}}(\f_{d-k+1}|\gfi_{[d-k]})= \left|\frac{\partial}{\partial \f_{d-k+1}} F_{\f_{d-k+1}}(s_k|\gs_{[k-1]},\gs_{-[d]}; \gfi_{[d-k]})\right|
\een
is inspired by the definition of the FD for a real parameter. Some remarks useful in the sequel follow.

\noindent
i)
When $m=d=1$, so that an ancillary statistic is not needed, formulas  \eqref{fid-general} and \eqref{fid-general-2} reduce to
$h_s(\f)=\left|\partial F_{\f}(s)/\partial \f\right|$,
the original proposal of \citet{Fisher:1930}. \\
\noindent
ii)
 When $d>1$ but the parameter of interest is $\f_1$ only, it follows from \eqref{fid-general} that its FD is simply given by
\be
 h_{\gs}(\f_{1})= \left|\frac{\partial}{\partial \f_1} F_{\f_{1}}(s_d|\gs_{[d-1]},\gs_{-[d]})\right|,
 \ee
 which is based on the whole sample and is also a CD.
 A typical choice for  $S_{d}$ is given by the MLE $\widehat{\f}_1$ of $\f_1$ and thus, when $\widehat{\f}_1$ is not sufficient, one has to consider the distribution of $\widehat{\f}_1$  given the ancillary statistic $\gs_{-[d]}$ as done in Section \ref{sec:p-star}.

\noindent
iii)
The FD  in \eqref{fid-general} is generally not invariant under a reparameterization of the model unless
the transformation from $\gfi$ to $\gla=(\la_1, \ldots, \la_d)$ say, maintains the same increasing order of importance in the components of the two vectors and  $\la_k$ is a function of $\f_1, \ldots, \f_k$, for each $k=1,\dots,d$, i.e. $\gfi(\gla)$ is a lower triangular transformation.

In \cite{Veronese:2014}  it is shown that the univariate FD/CD for a real NEF is asymptotically normal. Because the   multivariate FD defined in  \eqref{fid-general}  is a product of one-dimensional conditional FDs, it is quite natural to expect that also the FD for a  $d$-dimensional NEF is asymptotically normal.
\begin{thm} \label{teo:asymptotic-normal}
Let $\gX=(\gX_1, \ldots, \gX_n)$ be an i.i.d. sample from a regular  NEF  on $\Rex^d$ with $\gX_i$ having density
$p_{\gte}(\gx_i)=\exp\{\sum_{k=1}^d \te_kx_k-M(\gte)\}$, mean vector $\gmu=\gmu(\gte)$ and variance function ${\bf V}(\gmu)=\mbox{Var}_{\gmu}(\gX_i)$. Furthermore, let $\bar{\gx}$ be  the observed value of the sample mean $\bar{\gX}=n^{-1}\sum_{i=1}^n \gX_i$.
If $\gX_i$ admits bounded density with respect to the Lebesgue measure or is supported by a  lattice, then the fiducial distribution of $\gmu$ is asymptotically order-invariant and asymptotically normal with mean $\bar{\gx}$ and covariance matrix ${\bf V(\bar{\gx}})/n$.
\end{thm}
Since $ \bf V(\bar{\gx})$ coincides with the reciprocal of both the observed and the estimated expected Fisher information matrix, recalling standard results about  asymptotic Bayesian posterior distributions, see e.g. \citet{johnson:1979}, the following corollary immediately holds.
\begin{cor}
Consider the statistical model specified in Theorem \ref{teo:asymptotic-normal}. If we assume a positive prior for $\gmu$ having continuous first partial derivatives, then  the asymptotic Bayesian posterior for $\gmu$ coincides with the asymptotically normal fiducial distribution.
\end{cor}
The asymptotic normality for multidimensional generalized fiducial distributions has been proved
by \citet{sonderegger:2014} under a set of regularity assumptions. We remark that the previous two results are specific for our definition of FD and  hold for NEFs without any extra regularity condition.  Furthermore, the proof of Theorem \ref{teo:asymptotic-normal}, given in the Appendix,  is completely different from the standard ones used to show asymptotic normality in  frequentist, Bayesian  or generalized fiducial settings. It is based on the convergence of the conditional distributions determined by the importance ordering of the parameters, it  heavily relies on the properties of the mixed parametrization of the NEF and consequently the result is given in terms of the mean parameter, which is more interpretable than the natural one.

Consider now a parameter  $\gla={\bf g}(\gmu)$, with ${\bf g}$ a one-to-one lower triangular continuously differentiable function.
From \citet[Prop. 1]{Veronese:2015}, it follows that the  FD for $\gla$ can be obtained from that for $\gmu$ by the standard change of variable technique and  thus we can construct the asymptotic FD in the same way.
However, Theorem \ref{teo:asymptotic-normal} states that the asymptotic FD for $\gmu$ is order invariant and hence it could be interesting to investigate if this is true also for an arbitrary parameter.  This conjecture might be reasonable looking at what happens in the Bayesian theory where the asymptotic (reference) posteriors do not depend on the order of the parameter components.
%
%
%
The following example illustrate this point.

{\em Example 4}.  Consider a sample of size $n$ from a multinomial experiment with outcome probability vector $\gp=(p_1, \dots, p_d)$, with $\sum_{k=1}^d p_k \leq 1$. Then, the vector of counts $\gS=(S_1, \dots, S_d)$, with   $\sum_{k=1}^d S_k \leq n$, is distributed according to a multinomial distribution with parameters  $n$ and $\gp$.
Using the step-by-step procedure described above, \citet[formula   25]{Veronese:2015} have proved that the FD for $\gp$
is a generalized Dirichlet distribution which depends on the specific fixed ordering of the $p_i$'s.
Assume now $d=2$ and consider the transformation  $\f_1=p_1/p_2$ and $\f_2=p_2$ which is not lower triangular. The FD of $\gfi=(\f_1,\f_2)$ in this order is
\ben \label{fid-es-4}
h_{\gs}(\gfi)\propto  \f_1^{s_1-1/2}(1+\f_1)^{-1/2} \f_2^{s_1+s_2-1/2}(1-(1+\f_1)\f_2)^{n-s_1-s_2-1/2}.
\een
This latter is different from  the FD induced by that of $\gp$ but  coincides with the posterior distribution obtained from the reference prior for $\gfi$, see \citet[Sec.   5.4]{Veronese:2015}.

Consider now the asymptotic setting. From  Theorem \ref{teo:asymptotic-normal} it follows that the asymptotic FD of $\gp= (p_1, \dots, p_d)$  is  N($\bar{\gx},  {\bf V(\bar{\gx}})/n)$, with $\bar{\gx}=\gs/n$ and  where the elements of $\gV(\bar{\gx})$ are $v_{kk}=\bar{x}_k (1- \bar x_k)$ and $v_{kr}=-\bar x_k \bar x_r$, $k\neq r$. It is easy to verify that for $d=2$ it induces on $\gfi$ a normal distribution with means $\bar x_1/\bar x_2$, $\bar x_2$, variances $\bar{x}_1(\bar{x}_1+\bar{x}_2)/(n\bar{x}_2^3)$, $\bar{x}_2 (1- \bar x_2)$ and covariance $-\bar{x}_1/\bar{x}_2$. This distribution coincides with the asymptotic distribution corresponding to \eqref{fid-es-4} ( derived for example using standard results on Bayesian theory) and this fact  supports our conjecture that asymptotic FDs are invariant to the importance ordering of the parameters and can always been derived through the standard delta method.

\section*{Appendix}

\noindent \emph{Proof of Theorem 1.}
For the sake of clearness, in this proof we denote by  $\Teh$ the MLE of a parameter $\te$   and by $\teh$ the corresponding estimate. If  $F_{\te}(\teh)$ is the distribution function of $\Teh$, assumed decreasing in $\te$, let  $1-F_{\te}(\teh)$ be the FD for $\te$.
If $F_{\te}(\teh)$ is increasing the proof is similar with $1-F_{\te}(\teh)$ replaced by $F_{\te}(\teh)$. Then
\ben \label{relaz-fiduciale}
H_{n,\teh}(z)=\mbox{Pr}_{\teh}\left \{\sqrt{n/b}(\te -\teh)\leq z \right \}=\mbox{Pr}_{\teh}\{\te \leq \te_n\}= 1- \mbox{Pr}_{\te_n}\{\Teh^*_n \leq \teh\},
\een
where $\te_n=z\sqrt{b/n}+\teh$ and  $\Teh^*_n$ is the MLE based on $n$ i.i.d. random variables $X^*_{n,i}, i=1,\ldots,n$, belonging to the same family of distributions of $X_i$, but with parameter $\te_n$.
Note that $\te_n$ converges to $\te$ for $n\rightarrow \infty$, because $\teh$ converges to the \lq\lq true\rq\rq{} value $\te$ for almost all sequences $(x_1,x_2,\ldots)$ and $\Te$ is an open interval. Thus  $\te_n$ belongs to $\Te$ for  $n$ large enough and for each $z \in \Rex$.
Starting from \eqref{relaz-fiduciale},  we can also write
\be
H_{n,\teh}(z)&=& 1- \mbox{Pr}_{\te_n}\{\sqrt{n/b}(\Teh^*_n-\te_n) \leq \sqrt{n/b}(\teh - \te_n)\}\\
&=& \mbox{Pr}_{\te_n}\{\sqrt{n/b}(\Teh^*_n-\te_n) \geq -z \}
=\mbox{Pr}_{\te_n}\{\sqrt{n/b}(\te_n-\Teh^*_n) \leq z\}.
\ee
Thus, the asymptotic expansion of $H_{n,\teh}(z)$ can  be derived by expanding the frequentist distribution function of $\sqrt{n/b}(\te_n-\Teh^*_n)$. This expansion can be directly obtained  by standard results, even if  $\{X^*_{n,i}, i=1,2,\ldots,n; n=1,2,\ldots,\}$ is a triangular array  because we consider only random variables and a first order approximation,
  see e.g. \citet{Garcia:1998} and  \citet[Theorem 5.22]{Petrov:1995}. The frequentist expansion of  $Z=\sqrt{n/b}(\te-\Teh)$ has been provided in several papers about matching priors under a set of regularity assumptions. Using  formula $(3.2.3)$ in \citet{Datta-Muk:2004} with $\te=\te_n$ and recalling that  $\Teh^*_n$ is the MLE of $\te_n$, we obtain
\ben \label{expansion-Z}
\mbox{Pr}_{\te_n}\{\sqrt{n/b}(\te_n-\Teh^*_n) \leq z\}=\Phi(z)-\phi(z)\left[
\frac{1}{2}\frac{I'(\te_n)}{I(\te_n)^{3/2}} + \frac{1 }{6} E_{\te_n}\left[\frac{\ell'''(\te_n)}{(-\ell''(\te_n))^{3/2}}\right] (z^2+2)\right]\frac{1}{\sqrt{n}} + O(\frac{1}{n}).
\een
Now,  because $-\lpp-I(\teh)= O_p(n^{-1/2})$  \citep[see e.g.][Sec. 3.5.3]{Severini:2000} and $\te_n-\teh=z\sqrt{b/n}=O_p(n^{-1/2})$,  we have  $I(\te_n)=-\ell''(\teh)+O_p(n^{-1/2})= 1/b+O_p(n^{-1/2})$.   Moreover, applying the delta method to the expectation in \eqref{expansion-Z}, this expansion becomes
 \be
\mbox{Pr}_{\te_n}\{\sqrt{n/b}(\te_n-\Teh^*_n) \leq z\}&=&\Phi(z)-\phi(z)\left[
-\frac{1}{2} b^{3/2} \lppp + \frac{1 }{6} b^{3/2}\lppp (z^2+2)\right]n^{-1/2} + O(n^{-1})\\
&=& \Phi(z)-\phi(z)\left[ \frac{1}{6}b^{3/2}\lppp(z^2-1)\right]n^{-1/2} + O(n^{-1}),
\ee
and the theorem is proved.  \st

\noindent
\emph{Proof of Corollary 1.}
The result follows immediately using the expansion of the posterior distribution provided by \citet[Theorem 2.1 and formulae (2.25) and (2.26)]{Johnson:1970}, assuming $\pi^J(\te)\propto I(\te)^{1/2}$ as prior. Notice that under the stated conditions on the posterior, this result can be used even if the prior is improper, as observed in  \citet[pag. 106]{Ghosh:2006}. \st

\noindent
\emph{Proof of Proposition 1.}  Recalling  that for a real NEF $\bar{x}=M^{\prime}(\teh)$, we can write
\be
 p^*_{\te}(\teh)=\exp\{n(\te M'(\teh)-M^*(\te))\},
\ee
 where $M^*(\te)=  \log(\int \exp\{n(\te M'(\teh))\} d \nu(\teh))$, with $\nu(\teh)$ denoting the dominating measure of the density of $\teh$.  Thus $p^*_{\te}(\teh)$ belongs to a regular real NEF and the result follows immediately by  \citet[Theorem 1]{Veronese:2014}.  \st

\noindent
\emph{Proof of Theorem 2.} Given a square $d \times d$ matrix $\bf A$, we use  ${\bf A}_{k[r]}$ to denote  the vector of the first $r$ elements of the $k$-th row of $\bf A$ and  ${\bf A}_{[k][k]}$ to denote the matrix identified by the first $k$ rows and columns of $\bf A$. Moreover, $ \bf A^T$ denotes the transpose of $\bf A$.and

In order to determine the asymptotic FD of $\gmu$ we apply the  step-by-step procedure introduced in Section 3 to the conditional distribution of $\bar{X}_{k}$ given $\bar{\gX}_{[k-1]}=\bar{\gx}_{[k-1]}$ for each $k$. Clearly for $k=1$, we have the marginal distribution of
$\bar{X}_{1}$. Since the covariance matrix $\gV(\gmu)$ of $\gX_i$ is finite, by the central limit theorem $\bar{\gX}$ is asymptotically N$({\gmu}, n^{-1} {\bf V(\bar{\gx})})$ and thus the marginal distribution of $\bar{\gX}_{[k]}$ is also asymptotically normal with $E(\bar{\gX}_{[k]})=\gmu_{[k]}$ and $Var(\bar{\gX}_{[k]})=n^{-1} {\bf V(\bar{\gx})_{[k][k]}}$.
Let
\begin{equation}\label{lambda_k}
\lambda_k=\mu_k+\gV({\bar{\gx}})_{k[k-1]} \left[\gV( {\bar{\gx}})_{[k-1][k-1]}\right]^{-1} (\bar{\gx}_{[k-1]}-\gmu_{[k-1]})
\end{equation}
and
\begin{equation}\label{qk}
q_{k}=\gV({{\bar{\gx}}})_{kk}-\gV( {{\bar{\gx}}})_{k[k-1]} \left[\gV({{\bar{\gx}}} )_{[k-1][k-1]} \right]^{-1} \left[\gV( {{\bar{\gx}}})_{k[k-1]}\right]^T.
\end{equation}
Using known results about the convergence of conditional distributions, see \citet[Theorem 2.4]{Steck:1957} or \citet[Sec.4]{Barndorff:1979}, it follows that
the conditional distribution of $\bar{X}_{k}$ given $\bar{\gX}_{[k-1]}=\bar{\gx}_{[k-1]}$   is asymptotically  N($\lambda_{k},n^{-1} q_{k}$).

Now recall that for a NEF it is always possible to consider the so called \lq\lq mixed parameterization\rq\rq{}   $(\gmu_{[k]},\gte_{-[k]})$ which is one-to-one with the natural parameter $\gte$, see e.g. \citet[ch. 3]{Brown:1986}.  For $\gte_{-[k]}$ fixed, the  distribution of $\bar{\gX}_{[k]}$ belongs to a NEF with parameter $\gte_{[k]}$ and thus the conditional distribution of $\bar{X}_{k}$ given $\bar{\gX}_{[k-1]}=\bar{\gx}_{[k-1]}$ depends only on $\theta_k$. The same must be  true of course for the corresponding asymptotic distribution, so that its mean parameter $\lambda_k$ depends only on $\te_k$ and hence only on $\mu_k$. Considering now the alternative mixed parameter $(\gmu_{[k-1]},\te_k,\gte_{-[k]})$,  it follows that there exists a one-to-one correspondence between $\te_k$ and $\mu_k$, for  $\gmu_{[k-1]}$ and $\gte_{-[k]}$ fixed.
As a consequence $\gmu_{[k-1]}$ can be fixed arbitrarily in the mixed parameterizations  $(\gmu_{[k-1]},\te_k,\gte_{-[k]})$ with no effect on the conditional distribution and we specifically assume  $\gmu_{[k-1]}=\bar{\gx}_{[k-1]}$. Using the parameter $(\bar{\gx}_{[k-1]},\mu_k,\gte_{-[k]})$, we have that $\la_k$ coincides with $\mu_k$, see \eqref{lambda_k}.
 Summing up,  each of the three parameters $\lambda_k$, $\te_k$ and $\mu_k$  represents a possible parameterization of the asymptotic conditional distribution of $\bar{X}_{k}$ given $\bar{\gX}_{[k-1]}=\bar{\gx}_{[k-1]}$, for fixed $\bar{\gx}_{[k-1]}$ and $\gte_{-[k]}$.
Thus we can find the asymptotic FD of $\lambda_k$. Consider now a random vector
$\bar{\gX}^*$ with distribution belonging to the same family of that of $\bar{\gX}$, with  mixed parameter  $(\bar{\gx}_{[k-1]}, \mu^*_k,\gte_{-[k]})$, where $\mu^*_k=\bar{x}_k+ z_k/\sqrt{n}$, with $z_k\in \Rex$, as in the proof of Theorem 1. Notice that the marginal distributions of $\bar{\gX}^*_{[k-1]}$ and of $\bar{\gX}_{[k-1]}$ are equal.
Such a $\mu^*_{k}$ is well defined  for  large $n$ since $(\bar{\gx}_{[k-1]}, \bar{x}_k,\gte_{-[k]})$
is a possible value for the mixed parameter in the distribution of the whole vector, because the NEF
is regular and thus the parameter space is open.

For   $n$ varying and fixed $k$,  the sequence of marginal sample means $\bar{\gX}^*_{[k]}$ derives from random vectors whose mean parameter depends on $n$, so that it forms  a triangular array.
In order to determine the FD of $\lambda_k$, we can consider the quantity   $\sqrt{n} (\lambda_k - \bar{x}_k)$, which is a sort of standardization of $\lambda_k$ in our fiducial context. Using (1), similarly to what done in \eqref{relaz-fiduciale}, we can write
\ben \label{inversione}
 \mbox{Pr}_{\bar{x}_k} \left( \sqrt{n} (\lambda_k - \bar{x}_k) \leq z_k | \bar{\gX^*}_{[k-1]} =  \bar{\gx}_{[k-1]}, \gte_{-[k]} \right) \hspace{-0.3cm}&=& \hspace{-0.3cm} \mbox{Pr}_{\bar{x}_k} \left( \lambda_k \leq  \bar{x}_k + \frac{z_k}{\sqrt{n}} \left| \bar{\gX^*}_{[k-1]} \right. =  \bar{\gx}_{[k-1]}, \gte_{-[k]} \right) \nonumber \\
  &=& \hspace{-0.3cm} 1-  \mbox{Pr}_{\la^*_k} \left( \bar{X}^*_k \leq  \bar{x}_k  | \bar{\gX^*}_{[k-1]} =  \bar{\gx}_{[k-1]}, \gte_{-[k]} \right) \hspace{-0.1cm},
\een
where $\la^*_k= \bar{x}_k + z_k/\sqrt{n}$.
Since $Var(\bar{\gX}^*_{[k]})$  is a continuous function of $\gmu^*=E(\bar{\gX}^*)$,  it converges to a positive definite matrix for each $k$ when  $\gmu^*$ converges to the \lq \lq true\rq\rq{} value of $\gmu$, for $n \rightarrow \infty$. Then, using the result on the convergence of a conditional distribution presented at the beginning of the proof with $\gmu$ replaced by $\gmu^*$, we have that $\bar{X}^*_{k}$ given $\bar{\gX}_{[k-1]}=\bar{\gx}_{[k-1]}$  is  asymptotically N($\lambda_k^*,q _{k}/n$).
Notice that from the existence of the second moment of each component of $\bar{\gX}^*_{[k-1]}$, it follows
that the condition required by \citet[Theorem 2.4, formula (28)]{Steck:1957}, for the case of triangular arrays, is satisfied.
Thus, the asymptotic normality of $\bar{X}^*_k$ given $\bar{\gX}_{[k-1]}=\bar{\gx}_{[k-1]}$ implies, for  $n \rightarrow +\infty$,
\be
\sup_{z_k} \left | \mbox{Pr}_{\la^*_k} \left( \bar{X}^*_k \leq  \bar{x}_k  | \bar{\gX}_{[k-1]}=\bar{\gx}_{[k-1]},\gte_{-[k]} \right)
     -\Phi\left (\sqrt{\frac{n}{q_k}}(\bar{x}_k -\la^*_k) \right )\right| \rightarrow 0 \quad a.s.
 \ee
  Recalling the expression of $\la^*_k$, we obtain
\be
\sup_{z_k} \left| \mbox{Pr}_{\bar{x}_k +z_k/\sqrt{n}}  \left( \bar{X}_k \leq \bar{x}_k | \bar{\gX}_{[k-1]} =\bar{\gx}_{[k-1]}, \gte_{-[k]}\right) - \Phi \left( -z_k/\sqrt{q_k}   \right) \right| \rightarrow 0  \quad a.s.
\ee
which, using \eqref{inversione}, gives
\be
 \sup_{z_k} \left|  \mbox{Pr}_{\bar{x}_k} \left( \sqrt{n} (\lambda_k - \bar{x}_k) \leq z_k | \bar{\gX}_{[k-1]} = \bar{\gx}_{[k-1]}, \gte_{-[k]}\right) - \Phi \left( z_k/\sqrt{q_k}   \right) \right| \rightarrow 0  \quad a.s.
\ee
We can conclude that the conditional FD of $\lambda_k $ given $\gte_{-[k]}$
is asymptotically normal with mean $\bar{x}_k$ and variance $n^{-1}q_k$, and thus it does not  depend on  $\gte_{-[k]}$. Recalling the one-to-one correspondence between $\theta_k$ and $\lambda_k$, for fixed  $\gte_{-[k]}$, and in particular that $\la_d$ is a one-to-one function of $\te_d$,   it follows that
 $\lambda_1,\lambda_2,\ldots,\lambda_d$ are asymptotically independent, so that the full vector $\gla=(\lambda_1,\lambda_2,\ldots,\lambda_d)$ is asymptotically N($\bar{\gx},n^{-1} \bf{Q}(\bar{\gx})$), where $\bf{Q}(\bar{\gx})$ is the diagonal matrix with $k$-th element $q_k$.

To obtain the asymptotic FD of $\gmu$ we  consider the one-to-one  lower-triangular transformation $\gmu=\textbf{g}(\gla)$, with $\gla=\textbf{g}^{-1}(\gmu)$ given by \eqref{lambda_k} for  $k=1, \ldots,d$.
Consider now the lower $d\times d$ triangular matrix $\bf A=\bf A(\bar{\gx})$ whose $k$-th row is made up by the  vector
 $\,-{\bf V(\bar{\gx})}_{k[k-1]} [{\bf V(\bar{\gx})}_{[k-1][k-1]}]^{-1}$, in the first $k-1$ positions, 1 in the $k$-th position and 0  elsewhere. Thus we can write  $\gla=\bf{A} \gmu + (\bf{I}-\bf{A})\bar{\gx}$ and $\gmu=\bf{A}^{-1} \gla +(\bf{I}-\bf{A}^{-1})\bar{\gx}$, with $\bf{I}$ denoting the identity matrix of order $d$.
 By applying the Cram\'er delta method it follows that $\gmu$ is asymptotically normal with (asymptotic) mean and covariance matrix $\bf{A}^{-1}\bar{\gx}+(\bf{I}-\bf{A}^{-1})\bar{\gx}=\bar{\gx}$ and $n^{-1} \bf{A}^{-1}\bf{Q}(\bar{\gx}) \bf{A}^{-1 \; T}$, respectively. We now show that $ \bf{A}^{-1}\bf{Q}(\bar{\gx}) \bf{A}^{-1 \; T}={\bf V(\bar{\gx}})$ or, equivalently, $\bf{Q}(\bar{\gx})=\bf{A}\bf{V}(\bar{\gx}) \bf{A}^T$.
 By direct computation it is easy to see that the $(k,h)$-th element of $\bf A V(\bar{\gx})$, $k,h=1,2,\ldots,d$, is
\begin{equation}\label{matrix AV}
{\bf V(\bar{\gx})}_{kh}-{\bf V(\bar{\gx})}_{k[k-1]} [{\bf V(\bar{\gx})}_{[k-1][k-1]}]^{-1}  {\bf V(\bar{\gx})}_{h[k-1]}^T.
\end{equation}
Notice that \eqref{matrix AV} is 0 for $k>h$ because the product of its last two factors gives a $(k-1)$-dimensional vector with 1 in the $h$-th position and 0  otherwise.
The  matrix $\bf A V(\bar{\gx}) A^T$ is of course symmetric, so that it is sufficient to proceed only for $k \geq h$. On its diagonal we have
\begin{equation}\label{matrix AVA^T}
{\bf V(\bar{\gx})}_{kk}-{\bf V(\bar{\gx})}_{k[k-1]} [{\bf V(\bar{\gx})}_{[k-1][k-1]}]^{-1}  {\bf V(\bar{\gx})}_{k[k-1]}^T, \quad k=1, \ldots,d,
\end{equation}
because the only nonzero element in the product of the $k$-th row of $\bf AV(\bar{\gx})$ and the $k$-th column of $\bf A^T$ is the product of  (\ref{matrix AV}), with $h=k$, and 1.  For $k>h$, the $(k,h)$-th element of $\bf AV(\bar{\gx})A^T$ is 0, because the first $k-1$ components of the $k$-th row of $\bf AV(\bar{\gx})$ and the last $d-h$ components of the $h$-th column of $\bf A^T$ are zero. Thus the matrix $\bf AV(\bar{\gx})A^T$ coincides with $\bf Q(\bar{\gx})$ and this completes the proof of the theorem.
\st

\section*{Acknowledgments}
This research was supported by grants from Bocconi University.

\bibliography{Asymp_fiducial_arXiv}

\end{document}